\documentclass{amsart}
\usepackage{amssymb}
\usepackage{amsthm}
\usepackage{comment}
\usepackage{hyperref}
\usepackage{amsmath}
\usepackage{amssymb}
\usepackage{enumerate}
\usepackage{graphicx}
\usepackage{amsthm}
\usepackage{amsopn}
\usepackage{amsfonts}
\usepackage{upgreek}
\usepackage{empheq}
\usepackage{caption}
\usepackage{xcolor}
\numberwithin{equation}{section}

\newtheorem{theorem}{Theorem}

\newtheorem{example}{Example}%
\newtheorem{remark}{Remark}%
\newtheorem{lemma}{Lemma}%
\newtheorem{definition}{Definition}%
\begin{document}

\title[A note on oscillatory behavior of DEPCAG and IDEPCAG]{A Note on the Oscillatory Behavior of Impulsive Differential Equations with Piecewise Constant Arguments via Difference Equations}

\author{Ricardo Torres Naranjo\textsuperscript{*}}
\address{Instituto de Ciencias Físicas y Matemáticas, Facultad de Ciencias, Universidad Austral de Chile\\
Campus Isla Teja s/n, Valdivia, Chile.}
\curraddr{}
\email{ricardo.torres@uach.cl}
\author{Eugenio Trucco Vera}
\address{Instituto de Ciencias Físicas y Matemáticas, Facultad de Ciencias, Universidad Austral de Chile\\
Campus Isla Teja s/n, Valdivia, Chile.}
\email{etrucco@uach.cl}
\author{\"{O}zkan \"{O}cal}
\address{Faculty of Science,
 Department of Mathematics, Akdeniz University, 07058
 Antalya, Turkey.}
\email{oocal@akdeniz.edu.tr}



\date{}

\dedicatory{Dedicated to the memory of Prof. Manuel Pinto Jim\'enez.}
\thanks{Corresponding author: Ricardo Torres Naranjo. Email: ricardo.torres@uach.cl}

\begin{abstract}
This paper studies the oscillatory behavior of solutions to linear nonautonomous impulsive differential equations with piecewise constant arguments, including both advanced and delayed cases
\[
x'(t) = a(t)x(t) + b(t)x([t-k]), \quad k \in \mathbb{Z}.
\]
By exploiting the hybrid structure of these systems, we reduce the problem to an associated difference equation whose coefficients explicitly incorporate both the continuous dynamics and the impulsive effects. 

Classical oscillation criteria for difference equations do not account for impulsive phenomena. Through the proposed reduction, we extend these criteria to a class of impulsive and non-impulsive equations (IDEPCA and DEPCA), obtaining explicit sufficient conditions for oscillation in terms of the original system data.

An example is provided to illustrate the applicability of the results.
\end{abstract}
\maketitle

\section{Introduction}

The study of differential equations with deviating arguments originates from the pioneering work of A.~Myshkis in the late 1970s~\cite{203}. This class of equations found one of its earliest applications in the work of S.~Busenberg and K.~Cooke~\cite{COOKE1}, who used such models to describe the dynamics of vertically transmitted diseases, illustrating their applicability in biological systems with delay-type effects.

As a general extension of Myshkis's ideas, M.~Akhmet introduced a broader family of functional differential equations of the form
\begin{equation}
    x^{\prime}(t) = g(t, x(t), x(\gamma(t))), \label{depcag_eq}
\end{equation}
where $\gamma(t)$ is called a \emph{piecewise constant argument of generalized type}.

The argument function $\gamma$ is defined in terms of two real sequences $\left(t_n\right)_{n \in \mathbb{Z}}$ and $\left(\zeta_n\right)_{n \in \mathbb{Z}}$ satisfying
$t_n \leq \zeta_n \leq t_{n+1},\,\, \text{for all } n \in \mathbb{Z},$
and
$\lim_{n \to \pm \infty} t_n = \pm \infty.$
Then,
\[
\gamma(t) = \zeta_n \quad \text{if } t \in I_n := [t_n, t_{n+1}),
\]
which means that $\gamma(t)$ is locally constant on each interval $I_n$.  
A simple but representative example is given by $\gamma(t) = [t]$, where $[\cdot]$ denotes the floor function, constant on each subinterval $[n, n+1)$ for $n \in \mathbb{Z}$.

Typically, each interval $I_n$ is divided into two subintervals corresponding to advanced and delayed arguments:
\[
I_n = I_n^+ \cup I_n^-, \qquad I_n^+ := [t_n, \zeta_n], \quad I_n^- := [\zeta_n, t_{n+1}].
\]
Equations of type~\eqref{depcag_eq} are known as \emph{Differential Equations with Piecewise Constant Argument of Generalized Type} (\emph{DEPCAG}).  
These systems are of hybrid nature, exhibiting both continuous and discrete dynamics. One of their most remarkable features is that solutions remain continuous even though the argument $\gamma$ is discontinuous.  
Under the assumption of continuity, it is possible to integrate \eqref{depcag_eq} over each interval $[t_n, t_{n+1})$ to obtain an associated difference equation, revealing a deep interplay between continuous and discrete dynamics (see~\cite{AK2, P2011, Wi93}).

If impulsive effects are introduced at the sequence of moments $\{t_n\}_{n \in \mathbb{Z}}$, one obtains an \emph{Impulsive Differential Equation with Piecewise Constant Argument of Generalized Type} (\emph{IDEPCAG}), which can be formulated as
\begin{equation*}
\begin{tabular}{ll}
$x^{\prime}(t) = g(t, x(t), x(\gamma(t))),$ & $t \neq t_n,$ \\[3pt]
$\Delta x(t_n) := x(t_n) - x(t_n^-) = H_n(x(t_n^-)),$ & $t = t_n, \quad n \in \mathbb{Z},$
\end{tabular}
\label{idepcag_gral}
\end{equation*}
where $x(t_n^-) = \displaystyle\lim_{t \to t_n^-} x(t)$ and each $H_n$ is a given impulsive operator (see~\cite{AK3, Samoilenko}).

When the argument $\gamma(t)$ is explicitly specified in the differential equation, the equation is referred to as a \emph{DEPCA} (Differential Equation with Piecewise Constant Argument).

\begin{definition}\hfill\\
A function $y(t)$ defined on $[\tau, \infty)$ is said to be \emph{oscillatory} if there exist two sequences $(a_n)$ and $(b_n)$ in $[\tau, \infty)$ such that $a_n \to \infty$ and $b_n \to \infty$ as $n \to \infty$, and
\[
y(a_n) \leq 0 \leq y(b_n), \quad \text{for all } n \geq M,
\]
for some sufficiently large integer $M$.  
In other words, $y(t)$ is oscillatory if it is neither eventually positive nor eventually negative. Otherwise, it is called \emph{non-oscillatory}.
\end{definition}

In the piecewise continuous framework, a function $y(t)$ may still be oscillatory even if $y(t) \neq 0$ for all $t \in [\tau, \infty)$. This is the perspective considered in the present paper (see~\cite{karakoc_berek}).

\begin{definition}{\cite{ELAIDY}}\hfill\\
A nontrivial solution $z(n)$ of a difference equation is said to be \emph{oscillatory} if, for every integer $M > 0$, there exists $n \ge M$ such that
\[
z(n)z(n+1) \le 0.
\]
Otherwise, $z(n)$ is called \emph{non-oscillatory}.
\end{definition}

\section{Aim of the work}
Using results from the theory of delayed and advanced difference equations, we derive sufficient conditions for the oscillation of solutions to the scalar linear nonautonomous \emph{IDEPCA}  
\begin{equation}
\begin{aligned}
x^{\prime}(t) &= a(t)x(t) + b(t)x([t- k]), && t \neq n, \,\, n\in\mathbb{N}, \\[0.5em]
\Delta x(n) &= c_n x(n^-), && t = n,
\end{aligned}
\label{sistema_idepcag_general_abstract}
\end{equation}  
where $t \geq 0$, $a(t)$ and $b(t)$ are locally integrable functions, $[\cdot]$ denotes the greatest integer function, and $k \in \mathbb{Z}$.\\

\noindent In this work, we propose oscillation criteria based on the discrete solution of \eqref{sistema_idepcag_general_abstract}, drawing on classical results from the theory of difference equations.\\
\noindent To the best of our knowledge, this theory has rarely been used to address key questions related to DEPCAG or IDEPCA equations. Because of the hybrid nature of these systems, difference equations offer a powerful yet largely unexplored toolbox for understanding their qualitative behavior.\\
\noindent We extend the results of \cite{TORRES_osc} by considering 
\[
\gamma(t) = n - k, \qquad t \in I_n = [n, n+1), \quad |k|\geq 1.
\]
\noindent Our work is organized as follows. We start with the basic material needed to understand the piecewise constant argument framework. We then review some key oscillation criteria from the theory of difference equations, which help us analyze oscillations in scalar linear DEPCA and IDEPCA equations. After that, we introduce new oscillation criteria for IDEPCA that extend classical results for \eqref{sistema_idepcag_general_abstract}. Finally, we include examples to illustrate how our results work in practice.

\section{Recent Advances on Oscillations in Differential Equations with Piecewise Constant Argument}

In their 1987 paper~\cite{Afta2}, A.~R.~Aftabizadeh, J.~Wiener, and J.-M.~Xu analyzed the linear delay differential equation with a piecewise constant argument given by
\begin{equation}
\begin{tabular}{ll}
$x^{\prime}(t) + a(t)x(t) + b(t)x([t-1]) = 0,$ & $t \geq 0,$
\end{tabular}
\label{depca_t_menos_1}
\end{equation}
where $a(t) \in C([0,\infty), \mathbb{R})$ and $b(t) \in C([0,\infty), \mathbb{R}^+)$, while $[\cdot]$ denotes the greatest integer function.  
The authors derived the following oscillation condition:

\begin{theorem}\label{t_menos_1_teo}\hfill\\
If
\begin{eqnarray*}
&&\limsup_{n \to \infty} \int_{n}^{n+1} \exp\!\left(\int_{n-1}^{s} a(u)\,du\right) b(s)\,ds > 1, \quad n \in \mathbb{N},
\end{eqnarray*}
then every solution of~\eqref{depca_t_menos_1} oscillates.
\end{theorem}

Later, in 1988, J.~Wiener and A.~R.~Aftabizadeh~\cite{WiAf88} examined a more general differential equation with a piecewise constant argument of the form
\begin{equation}
\begin{tabular}{ll}
$y^{\prime}(t) = a(t)y(t) + b(t)y\!\left(r\!\left[\frac{t+k}{r}\right]\!\right),$ & $y(0) = y_0,$
\end{tabular}
\label{idepca_afta_wiener}
\end{equation}
where $a(t)$ and $b(t)$ are continuous functions on $[0,\infty)$, and $k, r$ are positive integers satisfying $k < r$.  
They obtained the following oscillation criterion:

\begin{theorem}\label{afta_wiener}\hfill\\
If either of the following conditions holds:
\begin{eqnarray*}
&&\limsup_{n \to \infty} \int_{rn-k}^{rn} b(t)\exp\!\left(\int_{t}^{rn} a(s)\,ds\right)\!dt > 1,\\
&&\liminf_{n \to \infty} \int_{rn}^{r(n+1)-k} b(t)\exp\!\left(\int_{t}^{rn} a(s)\,ds\right)\!dt < -1, \quad n \in \mathbb{N},
\end{eqnarray*}
then every solution of~\eqref{idepca_afta_wiener} is oscillatory.
\end{theorem}

In 2013, K.-S.~Chiu and M.~Pinto~\cite{Kuo_pinto_2013} investigated the scalar differential equation with a generalized piecewise constant argument:
\[
y'(t) = a(t)y(t) + b(t)y(\gamma(t)), \quad y(\tau) = y_0,
\]
where \(a(t)\) and \(b(t)\) are continuous real-valued functions on \(\mathbb{R}\).  
Their analysis focused on the existence of oscillatory and periodic solutions.  
This work appears to be the first to systematically address such properties for equations involving generalized piecewise constant arguments.  
However, their study assumes the continuity of solutions and does not extend to impulsive settings such as IDEPCA or IDEPCAG.\\

In 2018, F.~Karakoc, A.~Unal, and H.~Bereketoglu~\cite{karakoc_berek} examined nonlinear impulsive differential equations with a piecewise constant argument (IDEPCA) described by
\[
\begin{cases}
y'(t) = -a(t)y(t) - y([t-1])f(z([t]) + h_1(y[t])), & t \neq n \in \mathbb{Z}^+,\\
z'(t) = -b(t)z(t) - z([t-1])g(z([t]) + h_2(z[t])), &
\end{cases}
\]
\[
\begin{cases}
y(n) = (1 + c_n)y(n^-),\\
z(n) = (1 + d_n)z(n^-),
\end{cases}
\quad n \in \mathbb{Z}^+,
\]
with initial conditions \(y(-1) = y_{-1},\; y(0) = y_0,\; z(-1) = z_{-1},\; z(0) = z_0\).  
Here \(a, b : [0, \infty) \to \mathbb{R}\) are continuous, \(f, g, h_1, h_2 \in C(\mathbb{R}, \mathbb{R})\), and \(c_n, d_n \in \mathbb{R}\) with \(c_n, d_n \neq 1\).  
The authors established both existence and uniqueness of solutions, and derived sufficient conditions for oscillation.\\

More recently, in~\cite{Kuo2023}, K.-S.~Chiu and I.~Berna (2023) analyzed impulsive differential equations with a piecewise constant argument of the form
\begin{equation}
\begin{tabular}{ll}
$y^{\prime}(t) = a(t)y(t) + b(t)y\!\left(p\!\left[\frac{t+l}{p}\right]\!\right),$ & $t \neq kp - l,$\\
$\Delta y(kp - l) = d_k y({kp - l}^{-}),$ & $t = kp - l,\; k \in \mathbb{Z},$
\end{tabular}
\label{kuo_2023}
\end{equation}
where \(y(\tau) = c_0\), \(a(t)\) and \(b(t)\) are continuous real-valued functions, \(p < l\), and \(d_k \ne -1\) for all \(k \in \mathbb{Z}\).  
They established oscillation criteria for the solutions of~\eqref{kuo_2023}.\\

In a more general framework, R.~Torres~\cite{TORRES_osc} (2024) studied the oscillatory behavior of nonautonomous impulsive linear differential equations with a generalized piecewise constant deviating argument:
\begin{equation}
\begin{aligned}
x^{\prime}(t) &= a(t)x(t) + b(t)x(\gamma(t)), && t \ne t_k,\\[0.3em]
x(t_n) &= (1 + c_n)x(t_n^{-}), && t = t_k,
\end{aligned}
\label{sistema_w_escalar}
\end{equation}
where \(x(\tau) = x_0\), \(1 + c_n \ne 0\) for all \(n \in \mathbb{N}\), and \(\gamma(t)\) denotes a generalized piecewise constant argument.  
The main oscillation criterion obtained there reads as follows:

\begin{theorem}\label{afta_torres}\hfill\\
Let \(k \ge M\) with \(M \in \mathbb{N}\) sufficiently large.  
If either of the following conditions holds:
\[
\text{if } 1 + c_{k(t)} > 0 \text{ and }
\begin{cases}
\displaystyle \lim_{k \to \infty} \sup_{k \in \mathbb{N}} \int_{t_{k(t)}}^{\zeta_{k(t)}} 
\exp\!\left(\int_{s}^{\zeta_{k(s)}} a(u)\,du\right)b(s)\,ds > 1,\\[0.5em]
\text{or}\\[0.5em]
\displaystyle \lim_{k \to \infty} \inf_{k \in \mathbb{N}} \int_{\zeta_{k(t)}}^{t_{k(t)+1}} 
\exp\!\left(\int_{s}^{\zeta_{k(s)}} a(u)\,du\right)b(s)\,ds < -1,
\end{cases}
\]
or
\[
\text{if } 1 + c_{k(t)} < 0 \text{ and }
\begin{cases}
\displaystyle \lim_{k \to \infty} \inf_{k \in \mathbb{N}} \int_{t_{k(t)}}^{\zeta_{k(t)}} 
\exp\!\left(\int_{s}^{\zeta_{k(s)}} a(u)\,du\right)b(s)\,ds < 1,\\[0.5em]
\text{or}\\[0.5em]
\displaystyle \lim_{k \to \infty} \sup_{k \in \mathbb{N}} \int_{\zeta_{k(t)}}^{t_{k(t)+1}} 
\exp\!\left(\int_{s}^{\zeta_{k(s)}} a(u)\,du\right)b(s)\,ds > -1,
\end{cases}
\]
then every solution \(x(t)\) of~\eqref{sistema_w_escalar} is oscillatory.
\end{theorem}

\begin{remark}\hfill\\
When the argument is given by $\gamma(t) = [t-k]$ with $k \in \mathbb{Z}$ and $|k| > 1$, we have
\[
\gamma(t) = n - k, \quad \text{for } t \in I_n = [n, n+1), \; n \in \mathbb{Z}.
\]
This case does not correspond to a generalized piecewise constant argument in the sense of~\cite{AK2, P2011}, since $\zeta_n = n - k \notin [n, n+1]$, but rather $\zeta_n \in I_{n-k} = [n-k, n-k+1)$.  
Hence, it was not considered in~\cite{TORRES_osc}, which dealt with oscillation criteria for the linear form associated with~\eqref{idepcag_gral}.  
However, the situations $\gamma(t) = [t]$ and $\gamma(t) = [t+1]$ were included in that study.

For the specific case $\gamma(t) = [t-1]$, Theorem~6 in~\cite{TORRES_osc} still applies, provided that $b(t) \ge 0$ in~\eqref{sistema_w_escalar}.  
In general, the configurations $\gamma(t) = [t]$ and $\gamma(t) = [t+1]$ are often more tractable, as they lead to first-order associated difference equations.
\end{remark}

\section{Definition of a solution for \texorpdfstring{\eqref{sistema_idepcag_general_abstract}}{(abstract system)} in the DEPCA and IDEPCA frameworks}

Let $\mathcal{P}_c(X, Y)$ denote the family of all functions defined from $X$ to $Y$ that are continuous for $t \neq t_k$ and left-continuous at the points $t = t_k$, where they may present first-kind discontinuities.  
Analogously, $\mathcal{P}_c^{1}(X,Y)$ designates the set of functions defined form $X$ to $Y$ whose derivative belongs to $\mathcal{P}_c(X,Y)$.

\begin{definition}[DEPCA solution of \eqref{sistema_idepcag_general_abstract}]\hfill\\
A function $x(t)$ is said to be a solution of the differential equation
\begin{equation}
    x'(t) = a(t)x(t) + b(t)x([t-k]), \quad k \in \mathbb{Z}, \label{depca_sistema_original}
\end{equation}
defined on the interval $[0,\infty)$, provided that the following conditions hold:
\begin{itemize}
    \item[(i)] $x(t)$ is continuous on $[0,\infty)$; 
    \item[(ii)] the derivative $x'(t)$ exists for all $t \in \mathbb{R}^+$ except possibly at discrete points $t = n$, $n \in \mathbb{N}$, where one-sided derivatives exist; 
    \item[(iii)] the equation \eqref{depca_sistema_original} is satisfied by $x(t)$ for all $t \in [n,n+1)$, and in particular, for the right-hand derivative of $x(t)$ at $t=n$, for every $n \in \mathbb{N}$.
\end{itemize}
\end{definition}

\begin{definition}[IDEPCA solution of \eqref{sistema_idepcag_general_abstract}]\hfill\\
A piecewise continuous function $y(t)$ is said to be a solution of \eqref{sistema_idepcag_general_abstract} if it satisfies the following properties:
\begin{itemize}
    \item[(i)] $y(t)$ is continuous on each interval $I_n = [n, n+1)$ and may have first-kind discontinuities at $t = n$, $n \in \mathbb{N}$. Moreover, $y'(t)$ exists for all $t \in [0,\infty)$ except possibly at $t = n$, where the one-sided derivatives exist, that is, $y(t) \in \mathcal{P}_c^{1}([n,n+1), \mathbb{R})$;
    \item[(ii)] on each interval $I_n$, $y(t)$ satisfies the ordinary differential equation
    \[
    y'(t) = a(t)y(t) + b(t)y(n-k);
    \]
    \item[(iii)] at the impulsive moments $t = n$, the jump condition
    \[
    \Delta y(n) = y(n) - y(n^-) = c_n(y(n^-))
    \]
    is fulfilled, that is,
    \[
    y(n) = (1 + c_n)y(n^-),
    \]
    where $y(n^-)$ represents the left-hand limit of $y(t)$ at $t=n$, for every $n \in \mathbb{N}$.
\end{itemize}
\end{definition}

\section{Existence and Uniqueness of Solutions for \texorpdfstring{\eqref{sistema_idepcag_general_abstract}}{Eq.}}

In this section, we establish the existence and uniqueness of solutions to the nonautonomous homogeneous linear IDEPCA \eqref{sistema_idepcag_general_abstract}:
\begin{equation*}
\begin{tabular}{ll}
$z^{\prime}(t) = a(t)z(t) + b(t)z([t-k]),$ & $t \neq n,\, n \in \mathbb{N},$ \\ 
$\Delta z|_{t=n} = c_n z(n^{-}),$ & $t = n,$ \\
\end{tabular}\end{equation*}
where \( k \in \mathbb{Z} \), \( z(t) \in \mathbb{R} \), \( t \geq 0 \), and the functions \( a(t) \) and \( b(t) \) are locally integrable. The sequence \( (c_n)_{n \in \mathbb{N}} \) is real-valued with \( 1 + c_n \neq 0 \) for all \( n \in \mathbb{N} \). \\

Since the notion of oscillation involves asymptotic behavior, we are only concerned with the forward continuation of solutions. \\

To prove existence and uniqueness, we proceed by explicit construction.

\begin{lemma}[Existence and Uniqueness]\hfill\\
There exists a unique solution of \eqref{sistema_idepcag_general_abstract} defined on \( [0, \infty) \).
\end{lemma}

\begin{proof}
Let
\[
\phi(t,s) = \exp\left(\int_{s}^{t} a(u)\,du\right)
\]
be the fundamental solution of the ordinary differential equation \( x^{\prime}(t) = a(t)x(t) \).

Let \( t \in I_n = [n, n+1) \) for some \( n \in \mathbb{N} \). On this interval, equation \eqref{sistema_idepcag_general_abstract} takes the form
\[
z^{\prime}(t) = a(t)z(t) + b(t)z(n-k).
\]
Hence, for \( t \in [n, n+1) \), the unique solution is given by
\begin{equation}
z(t) = \phi(t, n)z(n) + z(n-k) \int_{n}^{t} \phi(t, s)b(s)\,ds.
\label{variacion_parametros_general}
\end{equation}

Evaluating at \( t = n+1 \), we obtain
\begin{equation}
z((n+1)^-) = \phi(n+1, n)z(n) + z(n-k) \int_{n}^{n+1} \phi(n+1, s)b(s)\,ds.
\end{equation}

Applying the impulsive condition at \( t = n+1 \), we have
\[
z(n+1) = (1 + c_{n+1})z((n+1)^-),
\]
which yields
\begin{equation}
z(n+1) = (1 + c_{n+1})\left( \phi(n+1, n)z(n) + z(n-k) \int_{n}^{n+1} \phi(n+1, s)b(s)\,ds \right).
\label{apoyo_coeficientes}
\end{equation}

Therefore, we obtain the following difference equation:
\begin{equation}
z_{n+1} = a_n z_n + b_n z_{n-k},
\label{ec_diferencias_a_resolver_depca}
\end{equation}
where
\begin{align}
a_n &= (1 + c_{n+1}) \exp\left( \int_n^{n+1} a(s)\,ds \right), \label{an_y_bn}\\
b_n &= (1 + c_{n+1}) \int_n^{n+1} \exp\left( \int_s^{n+1} a(u)\,du \right) b(s)\,ds.\nonumber
\end{align}

For any admissible set of initial conditions, equation \eqref{ec_diferencias_a_resolver_depca} admits a unique recursively defined solution. Substituting the corresponding sequence \( (z_n) \) into expression \eqref{variacion_parametros_general}, we recover the unique solution of \eqref{sistema_idepcag_general_abstract}.
\end{proof}

\section{Relationship Between the Oscillatory Behavior of IDEPCA and Difference Equations}
Next, we present one of the main tools used in this work. Although elementary, it provides the fundamental idea for linking the oscillatory behavior of IDEPCA and DEPCA with the oscillation theory of difference equations:

\begin{lemma}\hfill\\
Consider the IDEPCA \eqref{sistema_idepcag_general_abstract}.  
If the discrete solution of \eqref{sistema_w_escalar} (that is, the sequence $\{z(n)\}_{n \in \mathbb{N}}$) is oscillatory, then the full solution $z(t)$ of \eqref{sistema_w_escalar} is also oscillatory.
\end{lemma}

\begin{proof}
Suppose, for the sake of contradiction, that the discrete sequence $\{z(n)\}_{n\in\mathbb{N}}$, given by
$$z_{n+1} = a_n z_n + b_n z_{n-k},$$
is oscillatory, but the full solution $z(t)$ is not.  

If $z(t)$ were non-oscillatory, then there would exist $T>0$ such that $z(t)$ remains nonzero and keeps the same sign for all $t \geq T$. In particular, for every integer $n \geq \lceil T \rceil$, we would have $z(n)$ sharing that same fixed sign and being nonzero. This would imply that the sequence $\{z(n)\}$ eventually keeps one sign, which contradicts the assumption that $\{z(n)\}$ is oscillatory. Therefore, $z(t)$ must also be oscillatory.
\end{proof}

\section{Some Useful Oscillatory Criteria from the Theory of Difference Equations}
\label{section_difference}

\noindent In this section, we introduce the main tools used to analyze the oscillatory behavior of solutions of the IDEPCA \eqref{sistema_idepcag_general_abstract}.  

\subsection{The difference equation \texorpdfstring{$z_{n+1} = a_n z_n + b_n z_{n-k}$}{zn+1 = an zn + bn zn-k} and its simplified form}

\noindent Consider the difference equation
\begin{equation}
z_{n+1} = a_n z_n + b_n z_{n-k}, \quad a_n \neq 0 \text{ for all } n \ge n_0, \quad k \in \mathbb{Z}, \; n \ge n_0.
\label{dif_equation_retardo_2}
\end{equation}

\noindent To simplify \eqref{dif_equation_retardo_2}, introduce a nonzero sequence $\alpha_n$ and define $y_n := \alpha_n z_n$, with $\alpha_{n_0} \neq 0$. Multiplying both sides by $\alpha_{n+1}$ gives
\[
y_{n+1} = \alpha_{n+1} a_n z_n + \alpha_{n+1} b_n z_{n-k}.
\]
Subtracting $y_n = \alpha_n z_n$, we get
$
y_{n+1} - y_n = (\alpha_{n+1} a_n - \alpha_n) z_n + \alpha_{n+1} b_n z_{n-k}.
$\\
\noindent Choose $\alpha_n$ so that the first term vanishes: $\alpha_{n+1} a_n - \alpha_n = 0 \Longrightarrow \alpha_{n+1} = \frac{\alpha_n}{a_n}$, hence
\begin{equation}\label{def_alpha}
\alpha_n = \alpha_{n_0} \prod_{j=n_0}^{n-1} \frac{1}{a_j}, \quad \alpha_{n_0} \neq 0.
\end{equation}

\noindent Substituting $z_{n-k} = y_{n-k} / \alpha_{n-k}$, we get $y_{n+1} - y_n = b_n (\alpha_{n+1}/\alpha_{n-k}) y_{n-k}$. Define
\begin{equation}\label{def_beta}
\beta_n := b_n \frac{\alpha_{n+1}}{\alpha_{n-k}}, \quad \Delta y_n := y_{n+1}-y_n.
\end{equation}
Then the reduced equation reads
\begin{equation}\label{ec_dif_reducida}
\Delta y_n + Q_n y_{n-k} = 0, \quad Q_n := -\beta_n, \; k \in \mathbb{Z}.
\end{equation}
In this way, $Q_n$ is defined as:
\[
Q_n =
\begin{cases}
-\,\dfrac{b_n}{a_n}, & k=0,\\[0.5em]
-\, b_n \displaystyle\prod_{j=n-k}^{\,n} \frac{1}{a_j}, & k>0,\\[0.5em]
-\, b_n \displaystyle\prod_{j=n+1}^{\,n+\ell-1}{a_j}, & k<0, \ \ell=-k.
\end{cases}
\]
\noindent Using the definitions of $a_n$ and $b_n$ given in \eqref{an_y_bn}, $Q_n$ can be expressed as
\begin{equation}\label{def_qn}
Q_n =
\begin{cases}
-\displaystyle 
\Big( \prod_{j=n-k}^{\,n-1} (1+c_{j+1}) \Big)^{-1} 
\int_n^{\,n+1} \exp\Big(-\int_{\,n-k}^s a(u)\,du \Big) b(s) \, ds, & k>0,\\
-\displaystyle 
\int_n^{\,n+1} \exp\Big(-\int_n^s a(u)\,du \Big) b(s)\, ds, & k=0,\\
-\displaystyle 
\Big( \prod_{j=n}^{\,n+\ell-1} (1+c_{j+1}) \Big) 
\int_n^{\,n+1} \exp\Big(\int_s^{\,n+\ell} a(u)\,du \Big) b(s)\, ds, & k<0, -k=\ell,
\end{cases}
\end{equation}

\subsection{Difference equations theory: some oscillatory criteria for \eqref{ec_dif_reducida}}\hfill\\
We now present several classical lemmata that provide essential oscillation criteria for equation \eqref{ec_dif_reducida},
\begin{align*}
\Delta y_n + Q_n y_{n-k}=0,\qquad k\in\mathbb{Z},
\end{align*}
which will serve as key tools in the sequel.
\begin{lemma}\label{lema_retardo_1}\cite{erbe}{(L.H.\ Erbe, B.G.\ Zhang, 1989)}\\
Consider equation \eqref{ec_dif_reducida} with $k\in\mathbb{N}$ and $Q_n>0$. If
\[
\liminf_{n\to\infty} Q_n = q > \frac{k^k}{(k+1)^{k+1}},\quad \text{ or } \quad \limsup_{n\to\infty} \sum_{j=n-k}^n Q_j  >1,
\]
then every solution of \eqref{ec_dif_reducida} is oscillatory. 
\end{lemma}

\begin{lemma}\label{lema_retardo_karpuz}\cite[Cor. 1]{Karpuz}{(B.\ Karpuz, 2017)}\\
Consider equation \eqref{ec_dif_reducida} with $k\in\mathbb{N}$ and $Q_n>0$. If
\[
\liminf_{n\to\infty} \left(\sqrt{Q_{n-1}}+\sqrt{Q_n}\right) >1,
\]
then every solution of \eqref{ec_dif_reducida} is oscillatory.
\end{lemma}

\begin{lemma}\label{lema_resumen}\cite[Th.~7.5.2]{Gyori_Ladas}{(I.\ Győri, G.\ Ladas, 1992)}\\
Consider the equation \eqref{ec_dif_reducida}, 
where $Q_n$ is a \textbf{sequence of real numbers} and $k \in \mathbb{Z}$.\\ Then, in each of the following cases, every solution of \eqref{ec_dif_reducida} is oscillatory:
\begin{enumerate}
    \item $k \geq 1$, $Q_n$ is eventually nonnegative, and
    \[
    \liminf_{n \to \infty} \sum_{i=n-k}^{n-1} Q_i >  \left( \frac{k}{k+1} \right)^{k+1}.
    \]
    \item $k = 0$ and $(1-Q_n)$ is not eventually positive.
    \item $k = -1$ and $(1+Q_n)$ is not eventually positive.
    \item $k \leq -2$, $Q_n$ is eventually nonpositive, and 
    \[
    \liminf_{n \to \infty} \sum_{i=n+k+1}^{n-1} (-Q_i) >  \left( \frac{k}{k+1} \right)^{k}.
    \]
\end{enumerate}
\end{lemma}

\begin{lemma}\label{lema_resumen2}\cite[Th.~7.5.3]{Gyori_Ladas}{(I.\ Győri, G.\ Ladas, 1992)}\\
Consider the equation \eqref{ec_dif_reducida} where $Q_n$ is a \textbf{sequence of real numbers} and $k \in \mathbb{Z}$.\\ Then, in each of the following cases, every solution of \eqref{ec_dif_reducida} is oscillatory:
\begin{enumerate}
    \item $k \geq 0$, $Q_n$ is eventually nonnegative, and
    \[
    \limsup_{n \to \infty} \sum_{i=n}^{n+k} Q_i > 1.
    \]
    \item $k \leq 1$, $Q_n$ is eventually nonpositive, and 
    \[
    \limsup_{n \to \infty} \sum_{i=n+k+1}^{n} (-Q_i) > 1.
    \]
\end{enumerate}
\end{lemma}

\begin{lemma}\cite{ocalan}\label{lema_avance2}{(Ö.\ Öcalan, Ö.\ Akin, 2007)}\\
Consider equation \eqref{ec_dif_reducida} with $k \leq -2$ and $Q_n\leq 0.$ If
\[
\limsup_{n \to \infty} Q_n=q < \frac{k^{k}}{(k+1)^{k+1}},
\]
then every solution of \eqref{ec_dif_reducida} is oscillatory.
\end{lemma}

\begin{remark}\hfill\\
We recommend the excellent book \cite{Gyori_Ladas} for the oscillation theory of delay difference and differential equations, the surveys \cite{Benekas},\cite{kikina} and \cite{stavroulakis} for readers interested in first-order delay difference equations, and the works \cite{chatzarakis} and \cite{ocalan} for those focusing on difference equations with advanced arguments.
\end{remark}

\section{Main Results}
Using Lemma 2 and the previously established discrete oscillation criteria, we deduce the oscillatory behavior of solutions to \eqref{sistema_idepcag_general_abstract}.
\begin{theorem}
Consider the following linear IDEPCA:
\begin{equation*}
\begin{aligned}
z'(t) &= a(t) z(t) + b(t) z([t-k]), & t \neq n,\\
z(n) &= (1+c_n) z(n^-), & t = n,
\end{aligned}
\end{equation*}
where $k \in \mathbb{Z}$, and
\begin{align*}
a_n &= (1 + c_{n+1}) \exp\Big(\int_n^{\,n+1} a(s)\,ds\Big),\\
b_n &= (1 + c_{n+1}) \int_n^{\,n+1} \exp\Big(\int_s^{\,n+1} a(u)\,du\Big) b(s)\,ds,
\end{align*}
with $Q_n$ and $\alpha_n$ defined in \eqref{def_alpha} and \eqref{def_qn}, respectively.\\
Then:
 \begin{itemize}
\item \label{lema_retardo_1_idepca}{(Erbe-Zhang IDEPCA oscillation criterion)}\\
Consider equation \eqref{ec_dif_reducida} with $k\in\mathbb{N}$ and $Q_n>0$. If
\[
\limsup_{n\to\infty}  \left(\displaystyle 
\Big( \prod_{j=n-k}^{\,n-1} (1+c_{j+1}) \Big)^{-1} 
\int_n^{\,n+1} \exp\Big(-\int_{\,n-k}^s a(u)\,du \Big) b(s) \, ds\right) < -\frac{k^k}{(k+1)^{k+1}},
\]
or
\[
\liminf_{n\to\infty} \sum_{j=n-k}^n \left( \left(\displaystyle 
\Big( \prod_{j=n-k}^{\,n-1} (1+c_{j+1}) \Big)^{-1} 
\int_n^{\,n+1} \exp\Big(-\int_{\,n-k}^s a(u)\,du \Big) b(s) \, ds\right)\right) <-1,
\]
then every solution of \eqref{ec_dif_reducida} is oscillatory. 
\item \label{lema_retardo_karpuz_2_idepca}{(Karpuz IDEPCA oscillation criterion)}\\
Consider equation \eqref{ec_dif_reducida} with $k\in\mathbb{N}$ and $Q_n>0$. If
{\small
\begin{align*}
&\sqrt{-\! \left(\displaystyle 
\Big( \prod_{j=n-k-1}^{\,n-2} (1+c_{j+1}) \Big)^{-1} 
\int_{n-1}^{\,n} \exp\Big(-\int_{\,n-k-1}^s a(u)\,du \Big) b(s) \, ds\right)} \\
&\quad +\sqrt{-\!\left(\displaystyle 
\Big( \prod_{j=n-k}^{\,n-1} (1+c_{j+1}) \Big)^{-1} 
\int_n^{\,n+1} \exp\Big(-\int_{\,n-k}^s a(u)\,du \Big) b(s) \, ds\right)} > 1,
   \quad \text{for all large } n,
\end{align*}
}
then every solution of \eqref{ec_dif_reducida} is oscillatory.
\item\label{lema_resumen_IDEPCA}{(Győri-Ladas IDEPCA oscillation criterion 1)}\\
Consider the equation \eqref{ec_dif_reducida}, 
where $Q_n$ is a \textbf{sequence of real numbers} and $k \in \mathbb{Z}$.\\ Then, in each of the following cases, every solution of \eqref{ec_dif_reducida} is oscillatory:
\begin{enumerate}
    \item $k \geq 1$, $Q_n$ is eventually nonnegative, and
    \[
    \limsup_{n \to \infty} \sum_{i=n-k}^{n-1} \left(\displaystyle 
\Big( \prod_{j=i-k}^{\,i-1} (1+c_{j+1}) \Big)^{-1} 
\int_i^{\,i+1} \exp\Big(-\int_{\,i-k}^s a(u)\,du \Big) b(s) \, ds\right) <-  \left( \frac{k}{k+1} \right)^{k+1}.
    \]
    \item $k = 0$ and $$1+\displaystyle 
\int_n^{\,n+1} \exp\Big(-\int_n^s a(u)\,du \Big) b(s)\, ds$$ is not eventually positive.
    \item $k = -1$ and $$1-(1 + c_{n+1}) \int_n^{\,n+1} \exp\Big(\int_s^{\,n+1} a(u)\,du\Big) b(s)\,ds$$ is not eventually positive.
    \item $k \leq -2$, $Q_n$ is eventually nonpositive, and 
    \[
    \limsup_{n \to \infty} \sum_{i=n+k+1}^{n-1} \left(\displaystyle 
\Big( \prod_{j=i}^{\,i-k-1} (1+c_{j+1}) \Big) 
\int_i^{\,i+1} \exp\Big(\int_s^{\,i-k} a(u)\,du \Big) b(s)\, ds\right) <  -\left( \frac{k}{k+1} \right)^{k}.
    \]
\end{enumerate}
\item \label{lema_resumen2_IDEPCA}{(Győri-Ladas IDEPCA oscillation criterion 2)}\\
Consider the equation \eqref{ec_dif_reducida} where $Q_n$ is a \textbf{sequence of real numbers} and $k \in \mathbb{Z}$.\\ Then, in each of the following cases, every solution of \eqref{ec_dif_reducida} is oscillatory:
\begin{enumerate}
    \item $k \geq 0$, $Q_n$ is eventually nonnegative, and
    \[
    \liminf_{n \to \infty} \sum_{i=n+k+1}^{n} \left(\displaystyle 
\Big( \prod_{j=i-k}^{\,i-1} (1+c_{j+1}) \Big)^{-1} 
\int_i^{\,i+1} \exp\Big(-\int_{\,i-k}^s a(u)\,du \Big) b(s) \, ds\right) < -1.
    \]
    \item $k \leq 1$, $Q_n$ is eventually nonpositive, and 
    \[
    \limsup_{n\to\infty}\sum_{i=n+k+1}^{n}(-Q_i)>1
    \]
    where
\begin{equation*}
Q_n =
\begin{cases}
-\displaystyle 
\Big( \prod_{j=n-k}^{\,n-1} (1+c_{j+1}) \Big)^{-1} 
\int_n^{\,n+1} \exp\Big(-\int_{\,n-k}^s a(u)\,du \Big) b(s) \, ds, & 0<k\leq1,\\
-\displaystyle 
\int_n^{\,n+1} \exp\Big(-\int_n^s a(u)\,du \Big) b(s)\, ds, & k=0,\\
-\displaystyle 
\Big( \prod_{j=n}^{\,n-k-1} (1+c_{j+1}) \Big) 
\int_n^{\,n+1} \exp\Big(\int_s^{\,n-k} a(u)\,du \Big) b(s)\, ds, & k<0.
\end{cases}
\end{equation*}

\end{enumerate}
\item \label{lema_avance_nooscilatorio_idepca}{(Öcalan-Akin IDEPCA oscillation criterion)}\\
Consider equation \eqref{ec_dif_reducida} with $k \leq -2$ and $Q_n\leq 0.$ If
\begin{align*}
\liminf_{n \to \infty} 
\,\,\left(\left(\prod_{j=n}^{\,n-k-1} (1 + c_{j+1})\right)
\left(\displaystyle \int_n^{n+1} \exp\left(\int_s^{n-k} a(u)\,du\right) b(s)\,ds\right)\right) >  -\frac{k^k}{(k+1)^{k+1}},
\end{align*}
then every solution of \eqref{ec_dif_reducida} is oscillatory.
\end{itemize}
\end{theorem}


\begin{proof}
By Lemma 2, it suffices to prove the oscillatory behavior of the difference equation \eqref{ec_dif_reducida}, since this implies the oscillatory behavior of \eqref{sistema_idepcag_general_abstract}. 

Now, consider \eqref{ec_diferencias_a_resolver_depca} and the quantities $\beta_n$ and $Q_n$ defined in \eqref{def_alpha} and \eqref{def_qn}, respectively. Observing that
\[
\limsup (-Q_n) = -\liminf Q_n, \qquad \liminf (-Q_n) = -\limsup Q_n,
\]
the result follows by applying Lemmas~\ref{lema_retardo_1}--\ref{lema_avance2}.
\end{proof}
\begin{remark}\hfill
\item If we set $c_n = 0, \forall n\in\mathbb{N}$, we obtain the DEPCA version of Theorem 4; that is, the non-impulsive cases of \eqref{sistema_idepcag_general_abstract}.
\end{remark}
\section{Some Examples of Oscillatory Linear DEPCA-IDEPCA Equations}

\begin{example}
Consider the following IDEPCA:
\begin{equation}
\begin{aligned}
x'(t) &= -\frac{1}{10} x(t) - t\,x([t-3]), && t \neq n, \\[0.3em]
x(n) &= \frac{8}{9} x(n^{-}), && t = n, \quad n \in \mathbb{N}.
\end{aligned}
\label{ejemplo_intro}
\end{equation}

The estimates of $a_n$ and $b_n$, as given by \eqref{an_y_bn}, are:
\begin{align*}
a_n &= \frac{8}{9} \exp\Big(-\int_n^{n+1} \frac{1}{10}\,ds \Big) = \frac{8}{9} e^{-1/10}, \\[0.2em]
b_n &= -\frac{8}{9} \int_n^{n+1} s \exp\Big(-\int_s^{n+1} \frac{1}{10}\,du \Big)\, ds\\[0.2em]
&= -\frac{80}{9} \bigl(1 - e^{-1/10}\bigr) 
\left(n + \frac{10 e^{-1/10} - 9}{1 - e^{-1/10}}\right) < 0, \quad \forall n \in \mathbb{N}.
\end{align*}

Hence, the corresponding discrete equation of \eqref{ejemplo_intro} is:
\begin{align*}
x_{n+1}
= \frac{8}{9} e^{-1/10} x_n -\frac{80}{9} \bigl(1 - e^{-1/10}\bigr) 
\left(n + \frac{10 e^{-1/10} - 9}{1 - e^{-1/10}}\right) x_{n-3}.
\end{align*}

By verifying the hypotheses of Lemma \ref{lema_retardo_1_idepca} (Erbe-Zhang IDEPCA oscillatory criterion), we compute:
\begin{align*}
\limsup_{n \to \infty} \Bigg( \left( \frac{9}{8} \right)^3 \Big( -\int_n^{n+1} s \exp\Big(-\int_{n-3}^s \frac{1}{10}\,du \Big) ds \Big) \Bigg) 
= -\infty.
\end{align*}

Therefore, every solution of \eqref{ejemplo_intro} is oscillatory.
\end{example}

\section*{Conclusions}
In this work, we develop a systematic way to reduce impulsive differential equations with piecewise constant arguments to associated difference equations, where the coefficients naturally reflect both the continuous dynamics and the impulsive effects.

This approach makes it possible to apply classical oscillation criteria from the theory of difference equations—originally formulated for non-impulsive discrete systems—to our setting. In particular, we obtain explicit conditions in terms of the original continuous and impulsive data, thereby extending these criteria to a class of impulsive differential equations where they are not directly applicable.

As a result, we provide a unified framework for studying oscillatory behavior in delayed, advanced, and impulsive cases, including situations that lie beyond the standard DEPCAG framework.
\section*{Disclosure statement}
The authors declare no conflict of interest. 
\section*{Acknowledgment}
Dedicated to the memory of Prof.\ Manuel Pinto Jiménez for his inspiring guidance. 

\bibliographystyle{plainurl}
\bibliography{biblio}
\end{document}